\documentclass[12pt]{amsart}
\usepackage{epsfig,epic,amssymb,eepic,color,graphicx}
\usepackage[all]{xy}
 \usepackage[T1]{fontenc}        
\newtheorem{thm}{Theorem}[section]
\newtheorem{proposition}[thm]{Proposition}
\newtheorem{prop}[thm]{Proposition}
\newtheorem{defn}[thm]{Definition}
\newtheorem{lemme}[]{Lemma}
\newtheorem{cor}[thm]{Corollary}

\newtheorem{remarque}[thm]{Remark}

\newtheorem{ex}[thm]{Example}
\newtheorem{rien}[thm]{}
 
\newcommand{\be}{\begin{enumerate}}
\newcommand{\ee}{\end{enumerate}}
\newcommand{\bi}{\begin{itemize}}
\newcommand{\ei}{\end{itemize}}

\def\R{\mathbb{R}}

\def\Z{\mathbb{Z}}

\def\om{\omega}
\def\Om{\Omega}

\def\ga{\gamma}    
\def\Ga{\Gamma}

\def\al{\alpha}
\def\be{\beta}

\def\vp{\varphi}

\def\la{\lambda}
\def\La{\Lambda}

\def\si{\sigma}
\def\Si{\Sigma}

\def\ep{\varepsilon} 
\def\mun{{^{-1}}}

\def\ds{\displaystyle}
\def\BB{{\mathcal B}}
\def\FF{{\mathcal F}}

\def\nd{\noindent}
\def\bull{${}$\hfill$\Box$\\}
\def\proof{\nd {\bf Proof.\ }}

\begin{document}
\vskip 1cm
\begin{center}
{\sc Haefliger's codimension-one singular foliations, \break
open books and twisted open books in dimension 3}
\end{center}

\title{}
\author{ Fran\c cois Laudenbach}
\address{Universit\'e de Nantes, UMR 6629 du CNRS, 44322 Nantes, France}
\email{francois.laudenbach@univ-nantes.fr}

\author{ Ga\"el Meigniez}
\address{Universit\'e de Bretagne-Sud,  L.M.A.M., 
BP 573, F-56017 Vannes, France}
\email{Gael.Meigniez@univ-ubs.fr}




\thanks{FL is supported by ANR Floer Power}

\begin{abstract} We consider  singular foliations of codimension one
on 3-manifolds, in the sense defined by Andr\'e Haefliger as being
 $\Gamma_1$-structures. We prove that under the obvious linear embedding 
condition, they are 
$\Gamma_1$-homotopic to a regular foliation carried by an open book
 or a twisted open book. 
The latter concept is introduced for this aim.
 Our result holds true in every regularity $C^r$, $r\geq 1$.
In particular, in dimension 3,
this gives
 a very simple proof of  Thurston's 1976 regularization theorem without using
Mather's homology equivalence.
\end{abstract}

\maketitle

\thispagestyle{empty}
\vskip 1cm

\section{Introduction}\label{section1}
\medskip

In this paper we prove a regularization theorem for
singular foliations of codimension one on closed 3-manifolds.
 They will even be 
homotopic  to some specific regular models.

\begin{thm} \label{reg}
Let $M$
be a closed 3-manifold equipped with a
 Haefliger structure $\xi$ of codimension one and
 differentiability class $C^r$, 
 $1\le r\le\infty$. Assume that the bundle normal to $\xi$ embeds
 into the bundle tangent to $M$.
Then $\xi$ is $C^r$-homotopic to a regular Haefliger structure
 (that is, a foliation) carried by an open book decomposition of $M$,
 or by a twisted open book decomposition of $M$.
\end{thm}

Let us recall from \cite{haef} that, on a manifold $M$, a Haefliger structure $\xi$
of codimension one, or  a $\Gamma_1$-structure, consists of:
\bi
\item  a line bundle $\nu=(E(\nu)\to M)$, that is, a real vector bundle of rank 1,
 called the bundle {\it normal} to $\xi$;
\item in the total space $E(\nu)$, a germ, still noted $\xi$, of codimension-one foliation along the zero section $i $ of $\nu$, 
 transverse to the fibers.
 \ei
The $\Ga_1$-structure $\xi$ is said to be {\it regular}
 when the foliation $\xi$ is also transverse to the
 zero section $i$. Then 
 the trace $\xi\cap i(M)$
 is a genuine foliation of $M$. If not, this trace is a singular foliation.

A \emph{homotopy} of $\xi$ is defined as a  
$\Gamma_1$-structure on $M\times [0,1]$ inducing $\xi$ on $M\times\{0\}$. A 
{\it regularization} of $\xi$
 is a homotopy
 to a regular $\Gamma_1$-structure. It does not exist in general. 
An obvious necessary condition is
that 
$\nu$ must embed into the tangent bundle $\tau M$. When $\nu$ is trivial and 
$\dim M=3$ this condition is fulfilled.

In what follows, the manifold $M$ will be  $C^\infty$ and  $\xi$ 
will be a $\Ga^r_1$-structure ($1\le r\le\infty$), meaning  that $\xi$
viewed as a foliation is
tangentially $C^\infty$ and transversely $C^r$,
that is, the foliation charts are 
$C^r$ in the direction transverse to the leaves.\\

Of course, ``carried by a (twisted) open book decomposition''
needs explanation. But we first comment on the
regularization aspect only.

Under this aspect,
our result is  a particular case of a general regularization 
theorem due to W. Thurston (see \cite{thurston}).
 Thurston's proof was based on a deep result due to J. Mather
 \cite{mather0,mather1}, namely,
the homology isomorphism between on the one hand the classifying space
 of the  group Diff$_c(\R)$ 
endowed with the discrete topology, and on the other hand
 the loop space $\Om$B$(\Ga_1)_+$.

We present here a proof of the above regularization theorem which does not
 depend on the latter result.
After the present research was achieved,
a regularization theorem in all dimensions, still
avoiding any difficult result,
was provided in \cite{meigniez}, without the open book models that
 we get here in dimension 3.  \\
 
Our proof
provides model foliations in all homotopy classes of $\Ga_1$-structure.
 In this introduction, to make short, we explain the model
 in the co-orientable
 case: that is, $\nu$ is trivial.

In this case,
the models are based on the notion of
 {\it open book decomposition}. 
 Recall that such a structure on $M$ consists of a link $B$ in $M$,
 called the {\it binding},
 and a fibration $p: M\smallsetminus B\to S^1$ such that, 
for every $\theta\in S^1$, the fiber
 $p^{-1}(\theta)$ is the interior of an embedded surface, 
called the {\it page} $P_\theta$, 
 whose boundary is the binding. The existence of an open book decomposition 
could have been
 proved by J. Alexander when $M$ is orientable, as a consequence
 of \cite{alex} (every orientable closed 3-manifold is a branched cover of 
the 3-sphere) and \cite{alex2} 
 (every link can be braided); but
 the concept was 
introduced by 
 H. Winkelnkemper in 1973 \cite{wink}. Henceforth,  we refer to the more 
flexible construction by E. Giroux, which includes the case
when $M$ is  non-orientable 
(see section \ref{open-books}).

 It is well-known that
every open book $\mathcal B$ gives rise to a foliation $\mathcal F_\mathcal B$,
 as follows.
The pages endow 
 $B$ with a normal framing.
 So, a tubular neighborhood $N(B)$ of $B$ is 
trivialized:
 $N(B)\cong B\times D^2$. Out of $N(B)$ the leaves of $\mathcal F_\mathcal B$
 are the pages modified by spiraling
 around $N(B)$; some neighborhood of $\partial N(B)$ in $N(B)$ is a union of compact leaves; and the 
rest of $N(B)$ is foliated by a Reeb component.

This foliation of $N(B)$ will be called a
\emph{thick Reeb component} and it is introduced
 for technical reasons in the homotopy argument of section \ref{homotopy}.
 This technical point could be avoided
 by using a theorem of F. Sergeraert \cite{serger} on Reeb components.
 We call such a foliation an {\it open book foliation}.
 
 The latter foliation can be modified by inserting a so-called 
{\it suspension foliation}. Precisely, let $K$ be
a compact subsurface in some leaf of $\mathcal F_\mathcal B$ out of $N(B)$, and 
let $K\times [-\ep,+\ep]$ be a {\it foliated} thickening  of it: each  
$K\times \{t\}$ is contained in 
  a leaf of 
 $\mathcal F_\mathcal B$. Let $\vp:\pi_1(K)\to {\rm Diff_c}(]-\ep,+\ep[)$
 be some 
 representation into the group of compactly supported diffeomorphisms;
$\vp$ is assumed to be trivial on the peripheral elements (necessarily
$\partial K $ is non-empty). It allows 
 us to construct a {\it suspension foliation} 
$\mathcal F_\vp$ on $K\times [-\ep,+\ep]$, whose leaves are transverse 
to the vertical segments $\{x\}\times [-\ep,+\ep]$ and whose holonomy
 is $\vp$. The modification consists of removing
 $\mathcal F_\mathcal B$ from the
 interior of
 $K\times [-\ep,+\ep]$ and replacing it by $\mathcal F_\vp$. The new foliation,
 denoted 
 $\mathcal F_{\mathcal B, \vp}$, is an {\it open book foliation modified by
 suspension.} 
 It is also
 said to be {\it carried} by the open book decomposition $\mathcal B$.\\
 
When working  with  a $\Ga_1$-structure $\xi $ whose normal bundle is twisted,
 it is necessary to introduce the notion of
{\it twisted open book}, which is inspired by work of E. Giroux
and will be explained in section \ref{open-books}. 
It is worth noticing that 
the existence of twisted open book decompositions on  $M$ with a given 
{\it normal bundle} $\nu$
follows from our regularization theorem applied to the so-called {\it twisted
trivial} $\Gamma_1$-structure with normal bundle
 $\nu$ (see definition \ref{trivial-twisted} and theorem \ref{twisted-thm}). \\

The paper is organized as follows (see the text for definitions).
 In section \ref{morsification},
 the given
$\Ga_1^r$-structure $\xi$ is modified by
 a \emph{$C^0$-small homotopy} to make it Morse (that is, with Morse singularities),
and to give it a \emph{pseudo-gradient} $X$ whose dynamics is trivial.
We deduce a decomposition of $M$ of the form $M=N(\Sigma)\cup H$,
where $N(\Sigma)$ is a tubular neighborhood of an embedded closed surface
$\Si$ which is transverse to $X$,
 and where $H$ is a handlebody with one or two connected components. Moreover,  $\Si$
is Poincar\'e-dual to $w_1(\nu)$, the first Stiefel-Whitney class of $\nu$. 

In section \ref{open-books},
we build an  open book ${\mathcal B}$ (possibly twisted) which is 
 somewhat adapted to $\Si$:
the  surface $\Si$ minus a disk is contained in a page,
 and the bundle normal to ${\mathcal B}$ is isomorphic to $\nu$. In passing, we obtain the
following theorem:
\begin{thm} \label{twisted-thm}Let $\nu$ be a line bundle over $M$
which embeds into the tangent space 
 $\tau M$. Then there exists a twisted (or non-twisted)
 open book decomposition $\mathcal B$ of $M$
 whose normal bundle is isomorphic to $\nu$. It is twisted whenever $\nu$ is so.
\end{thm}

In section \ref{plateau}, we apply to
 $\xi$ a new homotopy to put it into \emph{plateau form}, meaning that it is {\it trivial} in $H$, and
in the interior of
 $N(\Si)$,  it is transverse to $X$.

In section \ref{homotopy} a final homotopy puts $\xi$ into
the desired form $\mathcal F_{\mathcal B,\vp}$ carried by ${\mathcal B}$.\\

Many arguments of this paper become significatively simpler when 
$\nu$ is trivial.
 The reader who is mainly interested  in this case may   refer to
 \cite{Laudenbach-Meigniez}. Also, the reader interested in prescribing the homotopy class of the plane field tangent to the resulting foliation is refered to the same preprint.

We are very grateful to Vincent Colin,  \'Etienne Ghys and 
 Emmanuel  Giroux for  their comments, suggestions and explanations.\\

\section{Morsification of the singularities
 and dynamics of a pseudo-gradient}\label{morsification}
 
 For proving  theorem \ref{reg}, 
the setting is  a closed $3$-manifold $M$ endowed with 
 a $\Ga_1^r$-structure $\xi$, with $r\ge 1$. Nevertheless, the next proposition \ref{morse} 
 holds  true in any dimension.

Denote by $\nu$
 the normal bundle to $\xi$, by
 $E(\nu)$ its the total space, and by 
$i:M\to E(\nu)$ the zero section. Regard $\xi$ as a foliation on some
 neighborhood of $i(M)$ in $E(\nu)$.
 
 A point $x\in M$ is said to be
 a {\it singularity} of $\xi$ when the zero section $i$
 is not transverse 
to $\xi$ at $x$.
One says that  $\xi$  is {\it Morse}
when at each singularity the 
contact of $i$ with $\xi$ is quadratic
non-degenerate in some foliated chart where $i$ is $C^2$. In that case, the trace
$\xi\cap i(M)$ is a singular foliation which is locally defined by a Morse function.

For every
 section $s$ of $\nu$ of class $C^r$ which is $C^0$-close to $i$,
 the pullback $s^*\xi$ is
a $\Ga_1^r$-structure on $M$ homotopic to $\xi$. Such a homotopy of $\xi$ 
is said to be a 
 \emph{$C^0$-small homotopy}. 

\begin{proposition} After a $C^0$-small homotopy,
$\xi$ is Morse. 
\label{morse}
\end{proposition}

\nd {\bf Proof.} Of course, in case $r\ge 2$ this proposition is obvious:
any $C^r$-generic choice of $s$ in a small enough
 $C^0$-neighborhood of $i$ works.
In case $r=1$ one has the following method which works in any
 dimension of $M$.

 At the beginning the zero section is covered by finitely many
 boxes bi-foliated
with respect to $\xi$ and to the $\R$-fibers. These boxes form a $C^r$ atlas. 
In each chart,
the zero section becomes the graph of a real function of class $C^r$.
 We choose a 
triangulation $Tr$ of 
$M=i(M)$ so fine that each simplex lies in such a box.
Then, for each simplex $\si$ in $Tr$
we choose
a bi-foliated neighborhood $V_\si$ such that when $\tau$ is a face of
 $\si$ we have 
$V_\tau\subset V_\si$. We need the two following lemmas.\\

\begin{lemme}
 Let $D^k$ be the standard compact  $k$-disk and let $U$ be 
 an open neighborhood in $D^k$ of 
the boundary $\partial D^k$. Let $f:U\to\R$ be a function of class $C^1$ without 
critical points. Let $V$ and $W$ be compact collar neighborhoods
 of $\partial D^k$
in $D^k$ with:
$$W\subset int\, V\subset V\subset U$$
Then
there exists a real function $g: D^k\to\R$ such that:
\bi
\item $g$ coincides with $f$ on  $W$ and has no critical point in $V$;
\item $g$ is $C^\infty$ on the complement of $V$;
\item $g$ is a Morse function (which makes sense since $g$ 
has no critical points in $V$).
\ei
\end{lemme}

\proof Let $X$ be a continuous vector field on $U$ such that $X\cdot f>0$.
 Let $\al$ be  some smooth bump function on $D^k$ which equals 1 on $W$ and 
has  support  in $V$.
Let $c>0$ be less than the minimum of $X\cdot f$ on $V$.
Let $m>0$ be the maximum of $\vert X\cdot \al\vert$ on $V$.
We choose  a $C^\infty $ function
 $h:D^k\to \R$ whose restriction to $V$ is very close
to $f$ in the $C^1$ topology. More precisely, 
we require that the following inequalities hold on $V$:

\bi
\item $X\cdot h>c$
\item $m\vert f-h\vert<{\ds\frac c 2}$\,.
\ei 
The function $h$ which has no critical points on $V$ may be chosen to be
 a Morse function on $D^k$. Then
 $g:=\al f +(1-\al) h$ has the wanted properties.
Indeed, at every point of $V$,
 $$X\cdot g=\al X\cdot f +(1-\al)X\cdot h +(X\cdot \al)(f-h)$$
 which is larger than $c/2$.
 \bull

\begin{lemme} Let $D^k$ be a compact $k$-disk embedded in $M$. Let
$\tilde f$ be a real $C^1$ function
defined on some neighborhood $ U$ of $\partial D^k$
in $M$. 
The functions $\tilde f$ and $\tilde f\vert U\cap D^k$
are assumed to have no critical points. Let  $g:D^k\to\R$ be a Morse function
which extends $\tilde f\vert U\cap D^k$. 
Then, there are a neighborhood 
$V$ of $D^k$ in $M$ and a function $\tilde g: V\to \R$ with the following
 properties:
\bi
\item $\tilde g\vert D^k=g$; 
\item $\tilde g$ coincides with $\tilde f$ in a neighborhood of $\partial D^k$;
\item $\tilde g$ is a Morse function and its critical points are those of $g$. 
\ei
\end{lemme}

\proof Near each critical point of $g$, the extension $\tilde g$ is constructed
 by adding a non-degenerate quadratic form in the coordinates of a transversal
to $D^k$. Hence, the last condition is satisfied. Away from the critical points
we have to solve the following extension problem: we are given
 a submanifold $N$ and a function $g:N\to\R$ without critical points.
We have to extend $g$ to a neighborhood of  $N$,
 the germ of the extension being already given near $\partial N$.
 It is easy by using  the projection $p:V\to N$ of 
some   tubular neighborhood of $N$ in $M$, well chosen near $\partial N$,
and  setting $\tilde g=g\circ p$.
\bull

We now return to the proof of proposition \ref{morse}.
It is done by induction on the skeleta of $Tr$. Assume we already
 have a section $s$
with the following property: for each simplex $\tau$ of  dimension less than 
$k$, the restriction $(s^*\xi)\vert\tau$ is Morse,
and there is a neighborhood $U$ of the $(k-1)$-skeleton such that 
$s^*\xi$ is non-singular on $U\smallsetminus Tr^{[k-1]}$.
If $\si$ is a $k$-simplex, $s\vert\si$ viewed in  the foliated chart $V_\si$
is the graph of a function $f$ to which it is allowable to apply successively 
the two preceding lemmas. It is clear from the construction that the new section is $C^0$-close to the initial zero section $i$.\bull

For example, the Morsification applies to the
 {\it twisted trivial} $\Ga_1$-structure in the
 following sense.

\begin{defn} \label{trivial-twisted} Let $\nu$ be a line bundle over $M$, whose first
 Stiefel-Whitney class is viewed as 
 a representation $w_1(\nu):\pi_1(M)\to\Z/2\subset GL_1(\R)$. 
 The twisted trivial $\Ga_1$-structure $\xi_0$ with  normal bundle $\nu $
 is  the foliation of $E(\nu)$ which is the suspension of $w_1(\nu)$.
\label{trivial}
\end{defn}

Then, all points of $M$ are singular. In other words, 
the zero section is a leaf. For every loop 
$\ga$ in this leaf, its holonomy
 is non-trivial if and only if $\nu$ is a twisted bundle over $\ga$.\\


In what follows, $\xi$ is assumed to be Morse.
We are interested in the dynamics of 
a so-called {\it pseudo-gradient} of $\xi$, defined as follows.
A {\em twisted} vector field $X$ on $M$ is a $C^\infty$ section $X$
of $Hom(\nu, \tau M)$. Then,
the sign of $X\cdot\xi$ is well-defined at each point of $M=i(M)$ where $X$
 is transverse to $\xi$.

\begin{defn}
A pseudo-gradient for $\xi$ is a twisted vector field $X$ on $M$
 such that  the Lyapunov inequality $X\cdot\xi<0$ holds everywhere but
 the singularities.
\end{defn}

As $\nu $ is possibly non-orientable, one cannot distinguish
the index of a singularity from its co-index; 
but one can distinguish saddle points from  center points thanks to their 
 phase portraits.
>From each saddle point $s$ start two {\it separatrices}, that is, isolated 
orbits of $X$ that begin
at $s$.

\begin{proposition}\label{finite} 
After some 
$C^0$-small homotopy,
 $\xi$ is still Morse and admits a 
pseudo-gradient $X$ whose dynamics has no recurrence 
(that is, every orbit has a finite 
length) and 
no separatrix joins two saddle points.
\end{proposition} 

\proof  The existence of a smooth pseudo-gradient is easy
 to prove even if $\xi$ is $C^1$ only.
Indeed, near the singularities, there are smooth charts
 of the singular foliation and 
the usual negative gradient in Morse coordinates
 is convenient.
Away from them, the Lyapunov inequality allows 
one to approximate a $C^0$ gradient
by a smooth pseudo-gradient.
 Let $X_0$ be such a  pseudo-gradient; its dynamics 
is not controlled.

Finitely many mutually disjoint  open 
2-disks, $d_1,\ldots, d_N$, are chosen in regular leaves
of the trace $\xi\cap i(M)$ such that every orbit of $X_0$ crosses  at least one shrinked
disk $\frac 12 d_k$,  for some $k\in{1,\ldots,N}$. 
 Following  Wilson's idea \cite{wilson},   the zero section $i$ of
$\xi $ and $X_0$ are changed in a neighborhood $D^2\times [-1,+1]$ of
 each disk into a {\it plug} such that every orbit of the modified
 pseudo-gradient $X$
is trapped by one of the plugs. Here are a few more details.

In these neighborhoods  $D^2\times [-1,+1]$, whose last coordinate is $t$, the foliation 
$\xi\cap i(M)$ is defined by the function $t$.  By a  $C^0$-small homotopy
supported in $D^2\times [-1,0]$, we just change 
the zero section $i$ to a section $s$ 
such that the singular foliation $\xi\cap s(M)$ is made of a cancelling pair of singularities, center--saddle. More precisely, it is defined by a Morse function $f$ which coincides with $t$
near the boundary and which has  a  pair of critical points
of indices  $(2,3)$. In these cylinders, the new pseudo-gradient  $X$ is the gradient of $f$ for the flat
metric. 
By using $D^2\times [0,1]$ in a convenient way, it is easy to make 
 the plug have the mirror symmetry with respect to
$D^2\times \{0\}$: any orbit of $X$ entering the cylinder at $p\times \{-1\}$ is trapped
 if $\vert p\vert\leq \frac 34$ 
 or gets out the cylinder at $p\times \{+1\}$. In the cylinder $D^2\times [0,1]$, the singular foliation 
 $\xi\cap s(M)$ is defined by a Morse function with a pair of critical points of indices $(0,1)$.
 
Thanks to the  mirror symmetry, very much like in Wilson's paper, the new twisted 
vector field $X$
has a {\it finiteness property}:
\bi
\item[] {\it Each orbit of $X$ has a  finite  length, and  connects two singularities.}
\ei
Let us give a proof.  We consider any half orbit $\la$ of $X$, starting at some
 point $x$ outside the plugs. We have to prove that the $\om$-limit set $\om(\la)$ is
 a singularity of $X$.  Let
$\la_0$ be the half orbit of $X_0$ through  $x$, with the same germ as $\la$ at $x$.
If $\la$ is trapped by some plug, that is, if it ends at some singularity in this
 plug, then we are 
done. So, we may assume that $\la$ is not trapped by any plug. Due to the mirror
 symmetry, as
recalled before, whenever $\la$ enters a plug at $p\times\{\mp 1\}$, $p\in D^2$,
 it gets out at $p\times\{\pm 1\}$. In particular, outside the plugs, we have $\la=\la_0$.
 It remains to show that $\om(\la_0)$ is a singularity of $X_0$. Indeed, if not, 
 $\om(\la_0)$
 would contain a whole orbit of $X_0$, thus would intersect some disk $\frac 12 d_k$.
 Thus $\la$ would be trapped by the corresponding plug, a contradiction.

Unfortunately, the mirror symmetry 
 in the plugs  implies that $X$ has some separatrices connecting two saddle points. 
Nevertheless,  one can destroy these connecting separatrices by a perturbation of $X$. 
We claim that the finiteness property  is preserved when the perturbation is small 
enough. 

By contradiction, assume that $X$ is the $C^0$-limit of some sequence $(X_n)$ 
of twisted vector fields  such that each $X_n$ has an orbit segment 
$\la_n $ of 
 length $n$.
 Then the Hausdorff accumulation set of $(\la_n)$ in $M$ for $n\to+\infty$
 would contain 
either  a half orbit of $X$ of infinite length, which does not exist by the 
finiteness property; 
or a half infinite (or periodic) 
 broken orbit $\La$ of $X$: that is, $\La$ is an infinite (or periodic)
 sequence of orbits
connecting successive saddle singularities. 
Then, if   $\La$ enters a plug
at $p\times\{\mp 1\}$ it must get out at  $p\times\{\pm 1\}$. 
Just as in the above proof of the finiteness property for
$X$, we conclude that there is some  half infinite (or periodic)
 broken orbit $\La_0$ of $X_0$
such that $\La=\La_0$ outside the plugs.

Let $\la_0$ be one of the 
 orbits of $X_0$ which are contained in 
$\La_0$. Then $\la_0$ meets some disk  $\frac12 d_k$.
 Thus, $\La$ is trapped by the 
corresponding 
plug and ends at the corresponding center singularity. 
Since $\La$ is contained in the accumulation set of $(\la_n)$, 
infinitely many $\la_n$'s are trapped in the same way as $\la_0$ and their
lengths arer bounded, a contradiction.

\bull\\

 
>From the pseudo-gradient $X$ of $\xi$ given by proposition \ref{finite}
we deduce the following topological objects.
Let $G$ be the topological closure in $M$ of the
 separatrices 
of all saddle points. It is a graph whose vertices (resp.  midpoints of the edges)
are the center (resp. saddle) singularities of $\xi$.
Thus $G$ admits an arbitrarily small
 tubular neighborhood $H$ 
 whose boundary is transverse to $X$. Let $\hat M$  be the complement 
of $int\,H$ in $M$. Since each orbit of
 $X$ from a point of 
$\partial H=\partial \hat M $ has a finite length, 
it must  return to the boundary. Therefore, $\hat M$  is fibered 
over a surface $\Si$, $\rho: \hat M\to\Si$,
the fibers being intervals ($\cong [-1,1]$)
tangent to $X$. By taking a section we think of $\Si$ 
as a closed surface embedded
 in $M\smallsetminus H$ and $\hat M$ becomes a tubular neighbohood
 $N(\Si)$ of $\Si$ in $M$. By construction of $X$, the normal bundle 
to $\Si$ in $M$ is $\nu\vert\Si$.\\

\begin{proposition}${}$\label{graphe}

\nd {\rm 1)} The line bundle  $ \nu\vert G$ is orientable and, 
for a suitable choice of the orientation,
 $X$ enters $H$ along $\partial H$.
  
\nd {\rm 2)} When $\nu$ is orientable,
$G$ has two connected components and $\Si$ is two-sided.

\nd {\rm 3)} When $\nu$ is non-orientable, $G$ is connected and 
 $\Si$ is one-sided.
 
 \nd {\rm 4)}  If $\nu$ embeds into the tangent bundle $\tau M$, the Euler characteristic 
 of $H$ is even.
\end{proposition}

\proof 1) We orient each separatrix from its saddle end point to its
 center end point. Since the separatrices are transverse to $\xi$,
this is an orientation of
$\nu\vert G$ over the complement of the singularities.  
It is easily checked that this orientation
extends over the singularities. Thus $X$ becomes a usual vector field near 
$G$ transverse to $\partial H$.

\nd 2) When $\nu$ is a trivial bundle, by the isomorphism above-mentioned,
$\nu(\Si,M)$ is trivial and 
 $\Si$ is two-sided. Thus,
$\partial \hat M$ has two connected components and $G$ also does.

\nd 3) Assume $\Si$ be  two-sided.  Then $\hat M\cong [-1,1]\times\Si$ and  $G$ has
 two connected components since each connected component of $H$ has a 
connected boundary. In that case, $\nu $ is orientable since there are no arcs 
in $H$ from $\Si\times\{+1\}$ to $\Si\times\{-1\}$.

\nd 4) By assumption, the pseudo-gradient $X$ is homotopic to a non-vanishing one. 
Then the number of zeroes of $X$ is even.
\bull

As a conclusion of this section, we shall be able to  continue the proof of theorem \ref{reg}
with a $\xi$ having a pseudo-gradient $X$ for which  the following statement holds true.

\begin{cor} \label{conclusion-2}
The above decomposition $M=N(\Si)\cup H$ has the following properties:
\bi
\item[{\rm 1)}] The handlebody $H$ has two or one components depending on that  $\nu$  is trivial or not.
\item[{\rm 2)}] The restricted bundle $\nu\vert H$ is trivial; equivalently, $\Si$ is Poincar\'e dual
to $w_1(\nu)$. For a convenient  orientation, $X$ enters $H$ along $\partial H$.
\item[{\rm 3)}] The fibres of $N(\Si)\to \Si$ are contained in orbits of the pseudo-gradient $X$.
 
\ei
\end{cor}

\section{Open books, twisted open books
 }\label{open-books}
 E. Giroux explained to us a Morse theoretical approach to open book
 decomposition, which  is based on  \cite{giroux} and  is recalled below,
 up to a change of terminology. 
 An easy adaptation allows us to handle similarly the {twisted} open book decompositions.
 
 \begin{defn} Let $W$ be a compact 3-manifold  
 and $f:W\to \R$ be a Morse function which is constant and maximal
 on $\partial W$. 
 A  compact 
 surface $S$ properly embedded in $W$ (that is, $\partial S\subset \partial W$) will be 
 called a 
  Giroux surface for $f$ when $f\vert S$ is a Morse function having the same critical points and the same isolated local extrema as $f$. 
  \end{defn}

The saddle points of $f\vert S$ may have index 1 or 2 for $f$. The property for $f$ of having a 
Giroux surface is kept when deforming $f$ among the Morse functions, even when crossing the critical values; for instance, through a deformation ending at a {\it self-indexing} Morse function
 (the value of a critical point is its Morse index in $W$). This definition is made for $W\subset M$;
 of course, we are mainly interested
  in the case  when $W$ is closed ($W=M$).
  
Theorem III.2.7 in Giroux's article \cite{giroux}
states the following with different words: 

\begin{thm}Let $M$ be a closed 3-manifold (orientable or not). There exist a
self-indexing Morse function $f: M\to \R$ and a  surface $S$
which is a 2-sided Giroux
surface in $M$ for $f$.\label{giroux-thm}  
\end{thm}

In that case, $S$ is separating (see below) and this data
immediatly gives rise to an open book decomposition 
$\mathcal B$ in the sense that
is recalled in the introduction. Indeed, 
 let $N$ be the level set 
 $f^{-1}(3/2)$. The smooth curve $B:=N\cap S$ will be the binding of the 
open book
 decomposition we are looking for. It can be proved that the following holds 
 for every  regular value $a, \ 0<a\leq 3/2$\,:
 \begin{itemize}
 \item the level set $f^{-1}(a)$ is the union along their common boundaries 
of two surfaces, $N_1^a$ and $N_2^a$,
 each one being  
 diffeomorphic to the sub-level surface $S^a:= S\cap f^{-1}([0,a])$;
 \item the sub-level set
 $M^a:= f^{-1}([0,a])$ is divided by $S^a$ into 
two parts $P_1^a$ and $P_2^a$ which are isomorphic handlebodies (with corners);
 \item $S^a$ is isotopic to $N_i^a$ through $P_i^a$, for $i=1,2$, by an 
isotopy fixing
 its boundary curve $S^a\cap f^{-1}(a)$.
 \end{itemize}
 This claim  is obvious when $a$ is small 
 and the property is preserved when crossing the critical level 1. 
 In this way the handlebody
 $H_-:= f^{-1}([0,3/2])$  is divided by $S^{3/2}$ into two diffeomorphic 
 parts $P_i^{3/2}$,
 $i=1,2$, and we have $N= N_1^{3/2}\cup N_2^{3/2}$. 
 We take
 $S^{3/2}$, $N_1^{3/2}$ and $N_2^{3/2}$ as pages; they are isotopic relative to $B$
 through respectively $P_1^{3/2}$ and $P_2^{3/2}$. 
 In $H_+:=f^{-1}([3/2, 3])$, we do the same construction with
 the function $3-f$. The open book decomposition is now clear. \\
 
We now generalize the notion of  open book decomposition.

\begin{defn} \label{gen-opb}
A \emph{generalized} open book decomposition of the closed connected 3-manifold $M$
is a pair $\mathcal B=(B,\mathcal P)$ where
$B$ is a co-orientable link, the \emph{binding,}
and $\mathcal  P$ is a codimension-one foliation of $M\smallsetminus B$,
 whose leaves 
are called the \emph{pages},
satisfying the following properties:
\begin{enumerate}
\item The union of each page with the binding is a compact (topological) surface.
\item In a tubular neighborhood $N(B)\cong D^2\times B$ of the binding, there are cylindrical coordinates
$(r,\theta,\phi)$, where $(r,\theta)$ are polar coordinates in $D^2$ and
$\phi$ is the projection to $B$, and
$\mathcal  P\vert N(B)$ is the foliation defined by $\theta=const$.
\item There is at least one page of $\mathcal  P$
 whose intersection with $N(B)$ is a single
annulus $\theta=\theta_0$.
\end{enumerate}
\end{defn}

After properties (1) and (2) the space of pages $(M\smallsetminus B)/\mathcal P$ is a compact connected 1-manifold, that is: $S^1:=\R/2\pi\Z$ or $I: =[-1, +1]$. 
When the first case holds,
$\mathcal B$ is a classical open book decomposition. 

\begin{defn} \label{open-book-def}A generalized open book decomposition $\mathcal B$ of $M$ whose space of leaves
is the interval $I$ is called a \emph{twisted} open book decomposition.
\end{defn}

Given a generalized open book $\mathcal B$, consider the projection
$p:M\smallsetminus B\to(M\smallsetminus B)/\mathcal P$.
In the non-twisted case,
each meridian of the binding
is mapped by $p$ onto the space of pages $S^1$ as a regular cover. The property (3) above forces this
cover to be of degree one.
In the twisted case, property (3) above forces that
each meridian meets each regular leaf in two points only.

 So, in the twisted case, $p$ is a singular fibration (or \emph{Seifert fibration}) 
 from $M\smallsetminus B$
onto $[-1, +1]$ 
which has  two one-sided  exceptional surface
fibers $p^{-1}(\pm 1)$, and which is a proper smooth 
submersion over the open interval.
A non-exceptional page 
is a 2-fold covering of any exceptional page $p^{-1}(\pm 1)$; notice that this covering is trivial
over a collar neighborhood of $B$ in the considered exceptional page.
The union of an exceptional page with the binding $B$
 is a smooth surface with
 boundary. But the union of a non-exceptional page
with $B$ is a closed 
surface showing
(in general)
an angle along  $B$.

\begin{figure}
\input{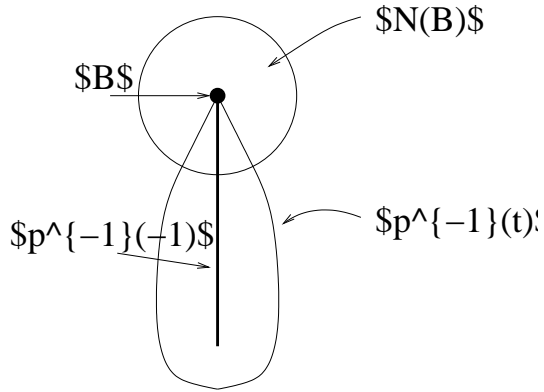}
\caption{Twisted open book.}
\label{book_fig}
\end{figure}


Notice that, since $B$ is co-orientable,
 a twisted open book $\mathcal B$
gives rise to a smooth foliation $\mathcal F_\mathcal B$
where each component of the binding is 
replaced with a Reeb component, the pages spiraling around it.

\begin{remarque}\label{external-hol}
{\rm
We may choose freely the {\it external holonomy} 
of this Reeb component,
that is, the germ by which the pages are rolled up around the binding, among the germs of
diffeomorphisms of $(\R,0)$ which are the identity on one side
and without fixed point on the other side. 
For the needs of the future homotopy argument at the end of our proof of theorem
 \ref{reg} (see section \ref{homotopy}), it is essential that this germ $\psi$ be chosen as a product of commutators. Of course, this assumption is known as being easily satisfied:
 for example, in $Aff(\R)$ the unit translation is a commutator; a classical conjugation 
  yields $\psi$ as above
 (see {\it e.g.} 
 \cite{meigniez-asian}, section 3). Moreover, 
 according to  Herman's theorem \cite{herman,herman2}, it is always satisfied;  but we want not to use this difficult result.}
 \end{remarque}

The {\it normal bundle}
to $\mathcal B$, noted $\nu(\mathcal B)$, is defined as being the normal 
bundle $\nu(\mathcal F_\mathcal B)$; it is well
defined up to isomorphism and is not thought of as a sub-bundle of $\tau M$,
 though  it embeds into 
the tangent bundle $\tau M$. Here is the simplest example of a twisted open book.\\

\begin{ex}{\rm Here, $M=S^1\times S^2$ is  thought of as the double of
 the solid torus
$S^1\times D^2$. The binding is the linear curve $B=(2,1)$ drawn on the flat  
 separating torus $T^2$.
One exceptional fiber is the 
 M\"obius band $\mathbb M$ in the solid torus, and its  boundary is $\partial \mathbb M= B$. 
The other exceptional fiber is its mirror copy. Notice that if $T^2$ is cut 
along $B$, one gets an annulus $A$ which doubly covers $\mathbb M$ when 
projecting along the normals to
$\mathbb M$. This covering is trivial over a collar neighborhood of $\partial \mathbb M$. }\\
\end{ex}

\begin{prop} \label{surf-to-ob}Let $f:M\to\R$ be a self-indexing (or ordered)
Morse function and $S$ 
be a one-sided Giroux surface for $f$. The following \emph{ middle condition} is assumed:
\bi
\item[{\rm (MC)}] The curve $B:=f^{-1}(\frac 32)\cap S$ is two-sided
 in its level set $f^{-1}(\frac 32)$.
\ei
Then there is a twisted open book $\mathcal B$ whose binding is the curve $B$ and 
such that $S\smallsetminus B$ is the union of  the two  exceptional pages. 
\end{prop}

\proof Observe first that $S^{3/2}:= S\cap f^{-1}([0, \frac 32])$ is one-sided. 
Indeed, if $a_0$ and $a_1$
are two points in a tubular neighborhood  $N(S^{3/2})$ which do not lie on $S$, they are joined
by an arc $\al$ in the complement of $S$ as $S$ is one-sided. By using the descending gradient 
lines of $f$, we can   push $\al$ into $N(S^{3/2})$ by an isotopy fixing the end points.

Knowing that the normal bundle to $S^{3/2}$ is twisted, it is possible to construct a piece
of twisted open book with binding $B$ inside $N(S^{3/2})$. The gradient of $f$
allows one to extend the open book structure on $f^{-1}([0, \frac 32])$ so that the complement of 
$B$ in
the level set $f^{-1}(\frac 32)$ is  a page. A similar construction
 in the upper part $f^{-1}([ \frac 32, 3])$ ends to build $\mathcal B$. \bull

 In general 
such a one-sided Giroux surface (or  a twisted open book) does not exist in $M$;
 the obstruction lies in the existence of a twisted line subbundle of
 $\tau M$ (compare theorem \ref{twisted-thm}).  Right now, we continue the proof of theorem
 \ref{reg}
 starting from the setting which has been stated in corolllary \ref{conclusion-2}. 
 
 We have a decomposition $M= N(\Si)\cup H$ where $\Si$ is a closed surface, $N(\Si)$
 is a $I$-bundle over $\Si$, whose projection is $\rho: N(\Si)\to \Si$, and $H$ is a handlebody with one or two components.
 Let $d_0$ be a small 2-disk in $\Si$; set $\Si'$ the closure of $\Si\smallsetminus d_0$
 and $M':= \rho^{-1}(\Si')$. 
 Consider the handlebody $H':=H\cup \rho^{-1}(d_0)$. We have a new decomposition
 $M=M'\cup H'$. We also recall the line bundle $\nu$ which is normal to the 
 $\Ga_1$-structure $\xi$ under consideration; by assumption it embeds into $\tau M$.\\
 
 \begin{prop}${}$\label{exist-ob}
 
 \nd {\rm 1)} When $\nu $ is trivial, there exists an open book $\mathcal B$
 such that $\Si'$ is contained in one page.
 
\nd{\rm 2)}  When $\nu $ is twisted, there exists a twisted  open book $\mathcal B$
 such that $\Si'$ is contained in an exceptional page and the normal  bundle $\nu(\mathcal B)$
 is isomorpic to $\nu$.\\
\end{prop}
 
 The first part is due to E. Giroux \cite{giroux2}.
 \smallskip

 \proof 
 For beginning with, assume that $\nu$ is twisted. The case when $\nu$ is trivial admits a similar treatment,
 with a few modifications which will be specified in the end.
According to proposition \ref{graphe} the handlebody $H$
is connected and its genus $g$ is odd; so, the genus of $H'$ is $g+1$ and even. A \emph{basis} of compression disks in $H'$
is a family of disjoint 
compression disks $\mathcal D =\{D_0, \ldots,D_{g}\}$ such that
  cutting H' along $\mathcal D$ gives rise to a 3-ball.
One passes from one basis to another by a sequence of elementary moves called \emph{slidings}.

Given such a basis, the disk $D_k$ is said to be {\it orientation-preserving} 
(resp. {\it orientation-reversing}) if cutting $H'$
along all the disks other than $D_k $ gives rise to a solid
torus (resp. a solid Klein bottle).

The following {\it sliding property} holds true: if  $\mathcal D$ is changed
 to $\mathcal D'$ by sliding $D_k$ over $D_j$ and if $D_k$ is orientation-reversing,
then the orientation type of  $D_j$ in $\mathcal D'$ is reversed and the other types remain unchanged.

\begin{lemme}\label{compression} 
The handlebody $H'$ admits a basis
$\mathcal D =\{D_0, \ldots,D_{g}\}$ of compression disks verifying:
\begin{enumerate}
\item[i)] $D_0$, \dots, $D_{g}$ are disjoint from 
the small disk $d_0=\Si\cap H'$;
\item[ii)] $D_1, \ldots,D_{g}$ are orientation-preserving;
 \item[iii)]  $H'$ splits into two connected domains $A_0$, $A_1$ whose
common boundary is either
$d_0\cup D_{0}$ (in case $H$ is orientable) or $d_0\cup D_{0}\cup D_{g}$ (in case $H$ is not orientable);
\item[iv)] For every even (resp. odd) $1\le k\le g-1$, the
disk $D_k$ is interior
to $A_0$ (resp. $A_1$);
\item[v)] In case $H$ is orientable, $D_g$ is interior to $A_1$.
\end{enumerate}
\end{lemme}

{\bf Proof} of lemma \ref{compression}.

i), ii) and iii). First
let $D_1$, \dots,
$D_g$ be any basis of compression disks for the handlebody $H$, and choose $D_{0}$
a compression disk for $H'$ parallel to $d_0$.

In case $H$ is orientable,  i)~--~iii) are immediate.

In case $H$ is non-orientable,
we shall change this basis by {sliding}
the disks one over another.
There is at least one orientation-reversing $D_k$, $k\ge 1$.
Sliding if necessary $D_k$ over $D_0$ and $D_{g}$,
one  makes $D_{0}$ and $D_g$ orientation-reversing.
Then,
sliding if necessary $D_{g}$ over $D_1$, \dots, $D_{g-1}$,
 one makes $D_1$, \dots, $D_{g-1}$ orientation-preserving.
Finally, sliding $D_0$ over $D_g$, one makes $D_g$ orientation-preserving.  Properties i)~--~iii) are verified.

iv) and v). One can pass any $D_k$, $1\le k\le g-1$, as well as $D_g$ in case $H$ is orientable, from $A_0$
to $A_1$, or from $A_1$ to $A_0$, by sliding twice $D_{0}$ over $D_k$. The orientation types are not changed. \bull\\

 We now start  constructing  a one-sided   Giroux surface satisfying (MC) and hence, according to
 proposition \ref{surf-to-ob}, a twisted open book.
 The surface $\Si'$, which is one-sided, is a Giroux surface in $M'$ with respect to some
 Morse function $f':M'\to \R$ having one minimum, $(g+1)$ critical points of index 1 and 
 which is constant on $\partial M'$. Now, we follow Giroux's algorithm for completing $\Si'$ to a closed Giroux surface. First, one has to attach 2-handles to $\partial M'$ in such a way 
 that: 
 \begin{enumerate} 
 \item[(vi)] each attaching curve intersects $\partial \Si'$ in two points exactly. 
 \end{enumerate}
 Thus, each 2-handle will
 produce simultaneously a 1-handle attached to $\Si'$, which allows one to extend both $f'$ 
 and the Giroux surface for it (\emph{cf.} \cite{giroux}). The previous disks $D_0, \ldots,D_g$ are devoted to be
  cores of these 2-handles, after convenient isotopies: in order that their attaching curves satisfy condition (vi), 
  pairs of intersection points with $\partial \Si'$ will be  created. 
  The following process is applied 
  in order to control the attachment of the last cell, after the $(g+1)$ surgeries of index 2.
  
 At the first step, a simple arc is drawn on $\partial M'\cap A_1$
from $\partial D_1$ to $\partial \Si'$, otherwise disjoint from the compression disks, and $\partial D_1$ is pushed by isotopy along this arc 
  in order to create two intersection points 
  with $\partial \Si'$. Let $(M_1,\Si_1)$ be the outcome of the  handle gluing to $(M',\Si')$.
   The boundary of 
  $\Si_1$ is made of an essential curve $c^1_0$ and a curve $c^1_1$ which bounds a disk in 
  $\partial M_1$.  By an isotopy of $M$ which preserves $H'$ and leaves 
  $D_0, D_2\ldots D_g$ fixed, one makes $c^1_0=\partial d_0$ and $c^1_1\subset A_1$.
  
\begin{figure}
\input{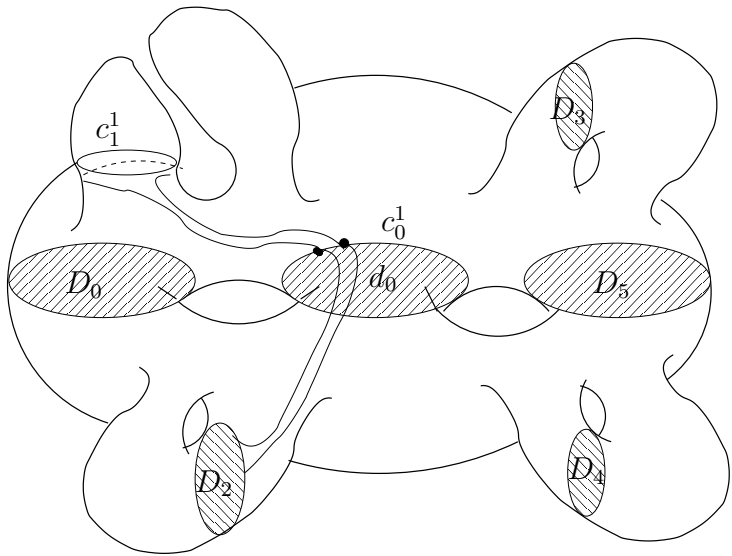}
\caption{The surface $\partial M_1$ and the isotopy
prior to the gluing of $D_2$, in case $g=5$
 and $H$ is not orientable.}
\label{sliding_fig}
\end{figure}

  At the second step, a simple arc is drawn on $\partial M_1$
from $\partial D_2$ to $c^1_1$, crossing 
 $c^1_0$ once, otherwise disjoint from the compression disks, and $\partial D_2$ is pushed by isotopy along this arc 
  in order to create two intersection points 
  with $c^1_0$ and also to surround
  $c^1_1$; condition (iv) guarantees that such an arc does exist (compare figure \ref{sliding_fig}).

  This surrounding amounts to a handle sliding of $D_2$ over $D_1$. According to condition 
  (ii), this operation does not change the type of the compression disks.
   And so on, until the gluing of $D_g$. 
   
  After step $g$, we have $(M_g,\Si_g)$ which is bounded by a torus or a Klein bottle
  whose complement in $M$ is a handlebody of genus 1
   for which $D_0$ is a compression disk. The boundary of $\Si_g$ is made
  of an essential curve $c_0^g=\partial d_0$, and a curve $c^g_1$ 
 which is a union of $g$ parallel circles bounding nested disks in $\partial M_g\cap A_1$.

 At the last step, a simple arc is drawn on $\partial M_g$
from $\partial D_0$ to $c^g_1$, crossing $c^g_0$,
 and $\partial D_0$ is pushed by isotopy along this arc
 in order to create two intersection points 
  with $c^g_0$ and also to surround $c^g_1$.  This amounts to a handle sliding of 
  $D_0$ over $D_1\cup \ldots\cup D_g$. By property ii) of lemma \ref{compression}, when $D_0$ is orientation reversing, all the other handles become orientation-reversing.
  After attaching this 2-handle, one gets 
  $M_{g+1}$ whose   complement in $M$ is a 3-ball  $\beta$,
   and a Giroux surface $\Si_{g+1}$
  in $M_{g+1}$
  whose boundary  is made of $(g+2)$ parallel circles. In order to close the Giroux surface by 
  only one  2-handle, it is necessary to have one circle only. We now explain this last step 
  of Giroux's algorithm.
  
  Let $\ga_0,\ga_1,\ldots, \ga_{g+1}$ be the boundary curves of the Giroux surface
  that we have in $\partial\beta$, the numbering being chosen so that two consecutive circles 
  bound an annulus avoiding the other circles.
   The regions of their complement in the 2-sphere are colored with two colors
  alternatively. A trivial 1-handle is attached to $M_{g+1}$ inside $\beta$ 
  whose core is a simple  unknotted arc $h_1$ in $\beta$
  having  one end point in $\ga_0$ and the other in $\ga_2$; moreover, one of the attached disks
   is turned by half a turn in order that the coloring extends along each side of
   the band which is attached to 
  $\Si_{g+1}$. Now, there is an obvious 2-handle which kills the previous 1-handle
  and whose core satisfies condition (vi) with respect to the previous Giroux surface.
  After this surgery, we have a Giroux surface in $M$ with a 3-ball removed, whose boundary
  is made of $g$ parallel circles in a 2-sphere. By repeating this operation we finally get the ideal 
situation where the boundary of the Giroux surface consists of one curve in the sphere
which bounds the {\it last cell} of  $M$: a 3-cell
 for closing $M$ containing 
  a 2-cell for closing the Giroux surface.

This construction yields a Giroux surface $S$ 
equipped with a Morse function $f$ whose critical values are not ordered with respect to 
their indices in $M$. We have added {\it trivial} 1-handles, $h_1, h_3,\ldots, h_g$ at a level
higher than critical points of index 2.

As in the classical Morse theory it is easy to make the reordering.  But, 
one has to take care of the {\it middle condition} (MC) in proposition \ref{surf-to-ob}.
 This condition could fail only when $H'$ is not orientable. We continue the proof in this case.
 
 The reordering mainly consists of extending the 1-handles $h_k$'s
 by the gradient lines of $f$ so that they are thought of as attached to $\partial \Si'$.
 Call $\bar h_1, \bar h_3,\ldots, \bar h_g$ these extended 1-handles.
 One checks that a case when (MC) certainly holds true is when they are compatible with
 the orientation $Or_0$ defined near $\partial \Si'$. However,
 by construction, they are compatible with the orientation $Or_\beta$ of $\beta$.
 As each compression disk of $H'$, after all the slidings we made, is orientation reversing,
 each time one crosses one $D_k$ when traversing $\partial \Si'$, the sign 
 $Or_\beta/Or_0$ changes. Fortunately, the shortest path in $\partial \Si'$  joining the feet
 of  $\bar h_1$ crosses exactly two compression disks,
 as the feet of $h_1$ are in $\ga_0$ and $\ga_2$ respectively. And similarly
 for the other $\bar h_k$'s. So, we are done.

Finally, having a Giroux surface and an ordered Morse function at hand, we have
 a twisted open book structure $\mathcal B$ with $ \Si'$ in an exceptional page.
  It remains to check that its normal bundle
  $\nu(\mathcal B)$ is isomorphic to $\nu$.
  
  The isomorphim class of a real line bundle is determined by its first Stiefel-Whitney class 
  $w_1(\nu)$, or by its Poincar\'e dual. By construction, $\Si$ is a Poincar\'e dual of $w_1(\nu)$
  with $\Z/2$ coefficients.  On the other hand, the Poincar\'e dual of 
  $w_1\bigl(\nu(\mathcal B)\bigr)$ is the Giroux surface $S$ that we have built. But, $d_0$
  and $S\smallsetminus int(\Si')$ are homologous in $(H',\partial H')$
  since they share a common boundary in the  aspherical manifold $H'$. Therefore, $S$ and $\Si$ are homologous and the considered bundles are isomorphic.\\

The case when $\nu$ is trivial is very similar. Here, $\Si$ is two-sided and the handlebody $H$
has two connected components, $H= H_1\sqcup H_2$, each one having the same 
 genus $g$. The handlebody 
$H'$ is made of the union of $H$ and the 1-handle $\rho^{-1}(d_0)$; it has genus $2g$.
Since the two handlebodies have diffeomorphic boundaries, if one is non-orientable, the other 
is neither.

There are compression disks $D_1, D_3,\ldots, D_{2g-1}$ in $H_1$ and 
compression disks $D_2, D_4,\ldots, D_{2g}$ in $H_2$. All together they are compression disks
 which cut $H'$ into one ball exactly. So, they will become the core of 2-handles attached 
to $M'$ after some suitable isotopy of their attaching curves. The boundary of $D_1$ is moved 
so that it intersects $\partial \Si'$ in two points. Next, we move the boundary of $D_2$
and so on. It is easy to check that the configuration of nested disks, similar to the previous 
case, can be realized at each step. \bull

\begin{cor} Theorem \ref{twisted-thm} holds true.
\end{cor}

\proof Start with the trivial or twisted trivial $\Ga_1$-structure $\xi_0$ according to 
whether its normal bundle is trivial or not (compare definition \ref{trivial-twisted}). Apply the Morsification process (section \ref{morsification}) until the decomposition 
$$M=N(\Si)\cup H$$
of corollary \ref{conclusion-2}. Thus, proposition \ref{exist-ob} yields the conclusion.\bull

\begin{rien}{\bf Modification by suspension.}\label{modif}
{\rm We end this section by  explaining how
the foliation $\FF_\BB$ associated above
to any twisted open book $\mathcal B$, 
can be modified by suspension in a similar manner as we did in section
\ref{section1} for open book foliations. 

Let $K$ be a compact subsurface in the exceptional
page $p^{-1}(-1)$. Consider a tubular neighborhood $N(K)$ of $K$ in $M$,
indeed a $[-\ep,\ep]$-bundle over $K$, which
 is compatible with $\mathcal F_\mathcal B$ in the following sense: the trace 
of $\mathcal F_\mathcal B$ on $N(K)$ is the 
foliation suspension of  the 
restricted  representation
$$\vp_0:=w_1(\nu):\pi_1(K)\to \Z/2={\rm Aut}([-\ep,\ep])
$$ (compare definition \ref{trivial}).
Let  $\vp: \pi_1(K)\to {\rm Diff}([-\ep,\ep])$ be  a representation
such that:
\bi 
\item[1)] for each $\al\in \pi_1(K)$ and  $x\in[-\ep,\ep] $ close
 to the end points, we have:
$$\vp(\al)(x)= \vp_0(\al)(x);$$
\item[2)] if $\al$ is peripheral, $\vp(\al)= \vp_0(\al)$.
\ei
These two  conditions allow us to remove $\mathcal F_\mathcal B$
from the interior of $N(K)$ and replace it with the suspension of $\vp$.
The modified foliation of $M$ is denoted by $\mathcal F_{\mathcal B,\vp}$. 
}\end{rien}

\begin{defn} 
 Any foliation of this form is said to be \emph{carried by the twisted open
book} $\mathcal B$.
\end{defn} 
After this definition, the statement of theorem \ref{reg} is meaningful.

\section{Homotopy to the plateau form}\label{plateau}

This section is a more step toward proving the regularization theorem \ref{reg}.
We introduce the following definition.

\begin{defn}\label{complete}
 Given a 
 $\Gamma_1$-structure $\xi$ on a space $G$, by an
 {\it upper (resp. lower) completion}
 of $\xi$ one means a foliation $\mathcal F$
of $G\times(-\ep, 1]$ (resp. $G\times[-1,\ep)$),
 for some positive $\ep$, which is transverse to every fiber
 $\{x\}\times(-\ep,1]$ (resp. $\{x\}\times[-1,\ep)$), 
whose germ along
$G\times \{0\}$ is $\xi$, and such that $G\times \{t\}$ is a leaf of
 $\mathcal F$ for every $t$ close enough to $+1$ (resp. $-1$).\\

\end{defn}

\begin{prop} \label{lower} 
Every co-orientable  $\Gamma_1^r$-structure $\xi$ on a simplicial complex $G$ of 
dimension 1, $r\ge 1$, admits an upper (resp. lower)
completion of class $C^r$.
\end{prop}

\proof Let us show the lower completion.
 After a fine subdivision of the edges, one reduces oneself to the case where,
 over each edge $\al$ of $G$, the germ of foliation
 $\xi$ is given by the level sets of a real function 
$f(x,t)$, $(x,t)\in \al\times]-\ep,+\ep [$; this function is  smooth 
in $x$, is $C^r$
in $t$,  and satisfies $\frac{\partial f}{\partial t}>0$ everywhere.
 Moreover, we may assume 
that the completion is already given over the 0-skeleton of $G$.
At this point, we  argue as in proposition \ref{morse}:
keeping the germ of $f$ fixed along $\al\times\{0\}\cup G^{[0]}\times [-1,0]$
one can arrange that $f$ is $C^\infty$ near $t=-\ep$. Finally,
the dimension of $\alpha\times[-1,0]$ being $2$, we are reduce 
to build a smooth  line field which fulfills the statement; the flow lines of 
this line field are the leaves of the wanted foliation. This can be done
 by partition of unity. \bull
\medbreak
Now we come back to the proof of theorem \ref{reg}. Recall the Morse
$\Ga_1^r$-structure $\xi$ obtained in proposition \ref{finite}, with its normal bundle $\nu$, its pseudo-gradient $X$, and the associated decomposition $M=N(\Si)\cup H$ 
(proposition \ref{graphe}, corollary \ref{conclusion-2}).
\def\plat{{\rm plat}}
\begin{defn}\label{plateau-forme}
 A $\Ga_1$-structure $\xi_\plat$ on $M$ with normal bundle $\nu$
 is said to be in \emph{plateau form} with respect to this decomposition,
if:
\bi
\item $\xi_\plat$ is trivial over $H$;
\item $\xi_\plat$ is regular over the interior of $N(\Si)$ and,
in this domain, it is transverse 
 to the fibers of $N(\Si)\to\Si$.
\ei 
\end{defn}

In other words, $\xi_\plat$
  is a suspension foliation in $N(\Si)$ and trivial
in $H$.

\begin{prop}\label{prereg}
 The $\Ga_1^r$-structure $\xi$ is homotopic to one, noted $\xi_{\rm plat}$, in plateau form with respect to the decomposition  $M=N(\Si)\cup H$.
\end{prop}

\begin{figure}
\input{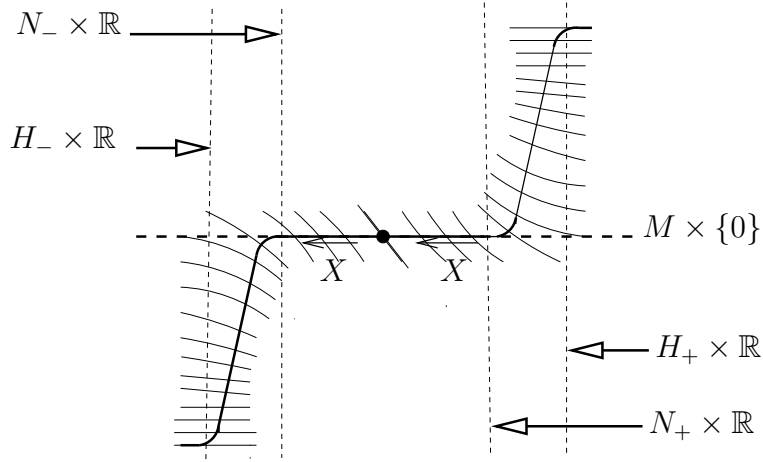}
\caption{
Homotopy to the plateau form tracked by some orbit of $X$.}
\label{plateau_fig}
\end{figure}

 This move is due to T. Tsuboi in \cite{tsuboi}, 
where it is given
as an exercise.
\medbreak
\proof   The homotopy of $\Ga_1$-structures will actually be a homotopy of the zero section in a foliated domain of the total space $E(\nu)$. We give the proof when $\nu$ is trivial
only, the other cases being similar. Then, $E(\nu)=M\times\R$, and $X$ is a  genuine vector field, and the graph $G$ 
formed by the separatrices of $X$
splits into $G_+\sqcup G_-$, where $G_+$ (resp. $G_-$) is a repeller (resp. an attractor)
 of $X$. By proposition \ref{lower}, 
 $\xi$ admits an upper (resp.  a lower) completion over $G_+$
(resp. $G_-$), and thus also over an open neighborhood $N_+$ (resp. $N_-$)
of $G_+$ (resp. $G_-$) in $M$. Recall that in section \ref{morse} $H$ was defined as an arbitrarily small handlebody neighborhood
of $G$ whose boundary is transverse to $X$. Thus, we can arrange
that the connected component $H_+$ (resp. $H_-$) of $H$ containing
$G_+$ (resp. $G_-$) is contained in $N_+$ (resp. $N_-$).
So we have a foliation $\mathcal F$ defined on a neighborhood of
$$
(M\times \{0\})\cup(H_-\times[-1,0])\cup(H_+\times[0,1])
$$ in $M\times\R$
which is transverse to $X$ on $\left(M\smallsetminus int(H_-\cup H_+)\right)\times\{0\}$ 
and tangent to  $H_\pm\times \{t\}$ 
for every $t$ close to $\pm 1$.

Recall (section \ref{morse}) that there is a diffeomorphism 
 $$F: M\smallsetminus Int(H'_-\cup H'_+)\to \Sigma\times[-1,+1]$$ which  maps orbit segments of $ X$ onto 
fibers.

For a small $\ep>0$, choose a function $\psi:{\R}\to[-1,+1]$
which is smooth, odd, and such that: 
\bi
\item $\psi(t)=0$ for $0\le t\le 1-3\ep$ and 
$\psi(1-2\ep)=\ep$; 
\item $\psi$ is affine on the interval $[1-2\ep,1-\ep]$; 
\item $\psi(1-\ep)=1-\ep$ and $\psi(t)=1$ for $t\ge 1$;
\item  $\psi'>0$ on the interval $]1-3\ep,1[$.
\ei

Let $s:M\to M\times\R$ be the graph of the function whose value is
 $\pm1$ on  $H_\pm$ and  $\psi(t) $ at the point 
$F^{-1}(x,t)$ for $ (x,t)\in \Sigma\times[-1,+1]$. When  $\ep$ is small
 enough,
it is easily checked that, for every $x\in \Si$, the path 
$t\mapsto s\circ F^{-1}(x,t)$ is transverse to $\mathcal F$ except at its end 
points.
Then,  $\xi_1:=s^*\mathcal F$ is homotopic to $\xi$ and obviously fulfills the
 conditions required in  proposition \ref{prereg}.
Indeed, at each point, $<s_*X,\xi>$ is a non-negative linear combination
of $<X,\xi>$ and $<-\frac{\partial\ }{\partial t}, \xi>$, hence, it is non-vanishing. (\emph{cf.} figure \ref{plateau_fig})\bull

 \section{Homotopy of $\Ga_1$-structures}\label{homotopy}

Now we complete the proof of theorem \ref{reg}.

On the one hand,
to the given $\Ga_1$-structure $\xi$ was associated
a decomposition $M=H\cup N(\Si)$, and according to
 proposition \ref{prereg}, $\xi$
 was homotoped to $\xi_{\rm plat}$ in plateau form with respect to this decomposition.
So,
$\xi_\plat$ is trivial on the
handlebody $H$, while on 
$N(\Si)$, it is the suspension foliation of some representation
$\vp: \pi_1(\Si)\to {\rm Diff}([-1,+1])$. The gluing with the trivial 
$\Ga_1$-structure on $H$ implies that, for each $\al\in\pi_1(\Si)$ 
and $x\in[-1,+1] $
near  the end points, we have:
$$\vp(\al)(x)=\vp_0(\al)(x),
$$where $\vp_0$ is obtained by restricting the first Stiefel-Whitney class 
of $\nu$ (compare condition 1) in \ref{modif}).
 Recall that
$d_0$ is a small disk in $\Si$, and $\Si'=\Si\smallsetminus int(d_0)$.

On the other hand, according to proposition \ref{exist-ob},  
$\Si'$ 
is contained in a
page of some open book $\mathcal B$; when $\nu$ is twisted,
$\mathcal B$ is twisted and $\Si'$ lies in an exceptional
page. 
Applying \ref{modif}, we modify 
$\FF_\BB$ by the suspension of the above-mentioned
 representation $\vp$  
and get a foliation $\FF_{\BB,\vp}$.
By construction, this foliation coincides with $\xi_\plat$ in $N(\Si')\cong \rho^{-1}(\Si')$,
where $\rho$ is the projection $N(\Si)\to\Si$.

Recall that $H'$ is the handlebody which is the 
union of $H$ and $\rho^{-1}(d_0)$. 
To establish theorem \ref{reg}, it remains to prove that 
the plateau $\Ga$-structure $\xi_{\rm plat}$ and the $\Ga$-structure $\FF_{\BB,\vp}$ carried by $\BB$ are homotopic as $\Ga_1$-structures on $H'$ relative to 
$\partial H'$. Notice that this statement is in fact independent of the representation 
$\vp$ since the modification of $\FF_\BB$ into $\FF_{\BB,\vp}$ did
 not modify the foliation in $H'$.

 The homotopy will be done in two steps,
 starting from $\FF_{\BB,\vp}$.

\subsection*{\bf First step.} 
First, one flattens its Reeb components. Set $N_1(B)$ the tubular neighborhood of $B$
which is the union of the Reeb components of $\FF_{\BB,\vp}$. As we are speaking here of 
non-thickened Reeb components we have $N_1(B)\subset N(B)$; the collar between both tubes is foliated by tori.

   \begin{lemme} \label{reeb}
There exists a homotopy, relative to $M\smallsetminus int(N_1(B))$,
 from $\FF_{\BB,\vp}$ to a $\Ga_1$-structure $\xi_\vp$  on $M$ which is trivial on $N_1(B)$.
 \end{lemme}
\proof 
Let $R$ be any connected component of $N_1(B)$.
Since the Reeb components of $\FF_{\BB,\vp}$ are \emph{thick}, the
holonomy of $\partial R$
 is trivial outside. This holonomy is generated by the
germ at $0$ of some self-diffeomorphism $\lambda$ of the real line,
whose support is contained in $[0,+\infty)$. Let $\FF_\lambda$ be the suspension of $\lambda$. It is 
a foliation on the annulus $S^1\times\R$, whose closed leaves are 
$S^1\times \{t\},\ t\leq 0$.

Then, on some small
neighborhood $N(R)$ of $R$, the foliation $\FF_{\BB,\vp}$ is the pullback
of $\FF_\lambda$ by some smooth map $F:N(R)\to S^1\times\R$.
More precisely, in cylindrical coordinates $(r, \theta,\phi)$
of $R\cong D^2\times S^1$ where
$\phi$ is the normal projection 
onto $S^1$ and $D^2$ is of radius 1, we take 
$F(r,\theta,\phi):=(\phi,\epsilon(1-r^2))$ with $\epsilon>0$ small.

Write $F(x)=(s(x),f(x))$ and define $G(x):=(s(x),g(f(x)))$ where $g$
 is any smooth function on the real line such that $g(t)=0$
 for every $t\ge 0$ and that $g'(t)>0$ for every $t<0$.
Then, $\xi_\vp:=G^*\FF_\lambda$ obviously works.
The homotopy
from  $\FF_{\BB,\vp}$ to $\xi_\vp$ is clear.
 \bull\\

 In the second step $\xi_\vp$ will be homotoped to $\xi_\plat$ 
relative to $\partial H'$.
 This second step will be different depending
 on whether  $\nu$ is trivial or not.
\subsection*{Second step, co-orientable case.}
Here the bundle $\nu$ is assumed to be trivial and $\mathcal B$ is a genuine open book decomposition.

  After the following lemma we
  shall be done with the homotopy problem.

 \begin{lemme} There exits a homotopy from $\xi_\vp$  to 
 $\xi_\plat$  relative to $N(\Si')$. 
 \end{lemme}

\begin{figure}[h]
\input{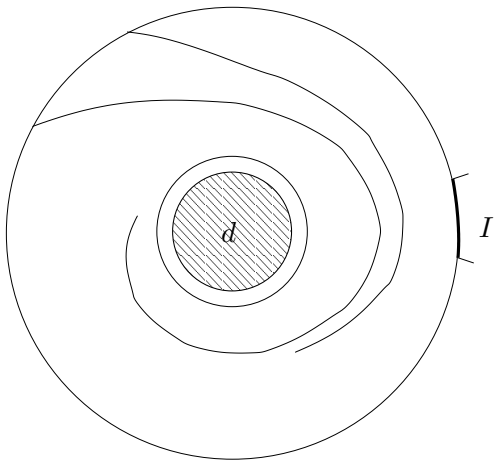}
\caption{The $\Ga_1$-structure $\xi_D$ on the 2-disk.}
\label{disk_fig}
\end{figure}

 \proof Recall the decomposition $M=N(\Si')\cup H'$.
 We have  to prove that the 
 restrictions of $\xi_\vp$ and $\xi_\plat$ to $H'$ are homotopic 
relative to $\partial H'$.
 Consider the standard closed 2-disk $D=D^2$ endowed with the $\Ga_1$-structure 
 $\xi_D$ which is shown on figure \ref{disk_fig}.
 
It is trivial on the small disk $d$ and regular on the annulus $D\smallsetminus
int(d)$. In the regular part,  the leaves are circles near $\partial d$ and the
 other leaves are spiraling, crossing $\partial D$ transversely. The
 restriction of $\xi_\vp$ to $H'$ 
has the form $f^*\xi_D$ for some map $f:H'\to D$.
Namely, on $N(B)\cong D\times S^1$ we take for $f$
the canonical projection onto $D$;
and we extend it continuously to $H'$ by sending every leaf
of $\xi_\vp\vert H'$, which is a subset of a page of $\BB$, to
a point of $\partial D$.
In particular, $f$ sends the 2-dimensional cylinder $\rho\mun(\partial d_0)$
to an interval $I$ embedded in $\partial D$.

In a similar way, $\xi_\plat\vert H'$ is the pullback of $\xi_D$
by the map $g:H'\to I \subset D$  which equals $f$ on $\partial H'$ 
 and sends $H$ to  $\partial I$
 (recall that $H$ has two connected components in the co-orientable case). Each regular leaf of $\xi_\plat$ is sent
 to the same point of $I$ as its trace on  $\rho\mun(d_0)$. Since $D$ 
 retracts by deformation onto $I$, consequently
 $\xi_\vp$  is homotopic to  $\xi_\plat$ on $H'$ relative
to $\partial H'$.

\bull\\

This finishes the proof of theorem \ref{reg} in the co-orientable case.

 \subsection*{Second step, twisted case.}

In the twisted case, the homotopy argument will be a little more sophisticated, because
 the classifying space of the Seifert fibrations over the interval is infinite-dimensional,
 as  the classifying space of the line bundles is. 

In the compact unit 3-ball $D^3$, let $\sigma$ be the orientation-preserving involution
of $D^3$ defined by 
$\si(x,y,t)= (x,-y,-t).$

\begin{lemme}\label{even}
There is a $\sigma$-invariant $\Ga_1$-structure $\xi_\sigma$ on $D^3$, whose restriction to the disk
 $(\partial D^3)\cap\{t\leq 0\}$ is conjugate to the
$\Ga_1$-structure $\xi_D$ represented on  figure \ref{disk_fig}.

\end{lemme}

\proof
It is more convenient to regard $D^3$ as the solid cylinder,
 namely, the product $D^2\times[-1,+1]$. In this model one has $\sigma(z,t)=(\bar z,-t)$, where $z=x+iy$.

The trace of $\xi_\sigma$
 on $\partial(D^2\times[-1,+1])$ 
will be trivial over $D^2\times\{\pm 1\}$
 and regular on $\partial D^2\times[-1,+1]$,
where it will indeed be the suspension of some diffeomorphism
$\zeta$ of the interval $[-1,+1]$. We first build this diffeomorphism.

Recall from  remark \ref{external-hol} that the rolling-up germ $\psi$
is chosen as a product of commutators. 
Let $0<\epsilon<1/2$. One easily makes a diffeomorphism $\gamma$ of the
interval $[-1,+1]$ such that:
\begin{enumerate}
\item[a)] One has $\gamma(t)\ge t$ with equality if and only if
 $t\ge 1-2\epsilon$ or $t\le -1+\epsilon$;
\item[b)] The germ of $\gamma$ at $-1+\epsilon$ is conjugate to $\psi$;
\item[c)] $\gamma$ is a product of commutators: $\gamma=[\alpha_1,\beta_1]\dots[\alpha_g,\beta_g]$. Here each $\alpha_j$, each $\beta_j$ is a diffeomorphism of the interval $[-1,+1]$ with support in $(-1,+1)$; and $[\alpha,\beta]$ denotes $\alpha\beta\alpha\mun\beta\mun$.
\end{enumerate} In order to construct such a $\gamma$ one starts  from a factorization
 of $\psi$ into $g$ commutators and extends each entry as a diffeomorphism of $[-1,+1]$.

  Then, set  $\zeta:=[\tau,\gamma\mun]$ where $\tau(t):=-t$. One has:
\begin{enumerate}
\item[i)] The germ of $\zeta$ at $-1+\epsilon$ is conjugate to $\psi$;

\item[ii)] $\tau\zeta\tau=\zeta\mun$. In other words, the suspension of $\zeta$ is $\sigma$-invariant in $\partial D^2\times[-1,+1]$;

\item[iii)] $\zeta(t)\ge t$. Indeed, for every
$u\in[-1,+1]$ one has $\gamma\mun(u)\le u\le\tau\gamma\mun\tau(u)$. Applying this at $u=\gamma(t)$ one has $t\le\gamma(t)\le[\tau,\gamma\mun](t)$.

\item[iv)] One has $\zeta(t)=t$ if and only if
$t=\gamma(t)$ and $\tau(t)=\gamma(\tau(t))$; that is, 
$t\le -1+\epsilon$ or $t\ge 1-\epsilon$.
\end{enumerate}

One will define $\xi_\sigma$ as the foliation given in $D^2\times[-1,+1]$ by the heigth function $h(z,t):=t$, modified as follows.
Create in the interior of $D^2\times[-1,+1]$ a number $g$ of pairs of singularities, the $j$-th pair ($j=1\dots, g$) consisting of two singularities
$s^j_1$, $s^j_2$ 
of respective indices $1$, $2$, in cancellation position. Let $f$ be the resulting Morse function.
One chooses its singular values so that  $f(s^j_1)<-1+\epsilon<1-\epsilon<f(s^j_2)$ ($j=1, \dots, g$). So, the intermediate level sets 
 have got some genus:  for $u\in [-1+\epsilon,1-\epsilon]$, the level set $f\mun(u)$ 
 is
a compact
surface of genus $g$ bounded by one circle. Its fundamental group 
$F_{2g} $ being non-abelian free on $2g$ generators, the $\al_j$'s and $\beta_j$'s define a representation $\la:F_{2g}\to {\rm Diff}([-1,+1])$. Next, 
in $f\mun[-1+\epsilon,1-\epsilon]$,  one changes the level surfaces of $f$ to the suspension of
the representation $\la$.
The resulting $\Ga_1$-structure on $D^2\times[-1,+1]$
induces the suspension of $\gamma$ on $(\partial D^2)\times[-1,+1]$. One pushes this structure by some convenient isotopy to make it coincide with the heigth function in the half cylinder $\{y\geq 0\}\times[-1,+1]$. 
Next, in the half cylinder  $\{y\leq 0\}\times[-1,+1]$, 
one performs the modification which is $\sigma$-symmetric to the preceding one. Finally, 
we  obtain  a $\sigma$-invariant $\Ga_1$-structure $\xi_\sigma$ in the solid cylinder.

Obviously $\xi_\sigma$ is trivial over $D^2\times\pm 1$, while its trace on $(\partial D^2)\times[-1,+1]$ is the suspension of the diffeomorphism $[\tau,\gamma\mun]$. Thus, thanks to property i) above,
the trace of $\xi_\sigma$ on the disk $(D^2\times-1)\cup(\partial D^2\times[-1,0])$ is conjugate to $\xi_D$.

\bull\\

The involution $\sigma$ is suspended to get a $D^3$-bundle
over $\R P^\infty$. Let $E$ be its total space.
By lemma \ref{even}, there exists a $\Ga_1$-structure $\xi_E$ on $E$
whose restriction to each fiber is $\xi_\sigma$. Its normal bundle
is the unique twisted bundle of rank one over $E$.

\begin{lemme} There are two continuous maps $f,g: H'\to E$, equal on $\partial H'$, such that
 $\xi_\vp\vert H'=f^*\xi_E$ and $\xi_\plat\vert H'=g^*\xi_E$.
\end{lemme}

\proof
Consider the circle $S^1:=\{t=0\}$ in $S^2$ in the coordinates considered above. As it is
 $\sigma$-invariant, it
defines a circle subbundle $E'\subset E$ over $\R P^\infty$.
Every circle fiber is transverse to $\xi_E$.
Actually $\xi_E\vert E'$
coincides with the (infinite-dimensional) horizontal foliation $\FF_{E'}$ on $E'$ which is the
suspension over $\R P^\infty$ of the orientation-reversing involution of the circle, $\sigma\vert S^1$.
The holonomy covering space of each leaf is $S^\infty$~, thus contractible.
So, $(E',\FF_{E'})$ is the  Haefliger classifying space for the 
groupoid generated by this involution
(\cite{haef2}). In other words, $({E'},\FF_{E'})$
 is the classifying space
for foliations of codimension one which are \emph{wandering}, that is, every leaf meets every transverse interval in a finite set.
Any such foliation $\FF$ on any manifold $V$ 
is the pullback of $\FF_{E'}$ by some \emph{classifying map} $V\to {E'}$~.
One also has the relative version: if $X\subset V$ is a submanifold
transverse to $\FF$
then every classifying map for $\FF\vert X$ extends to some classifying
map for $\FF$~.

Just as in the orientable case above, $\xi_\vp\vert N(B)$ is the pullback of $\xi_D$ through the canonical projection $N(B)\cong D^2\times S^1\to D^2$~. On the other hand,
in the complement $H'\setminus int(N(B))$ the $\Ga_1$-structure $\xi_\vp$ is a wandering foliation. Embed $D^2$ into $E$
 as one half of the boundary of some 3-ball
fibre, so that $\xi_D$ is the restriction of $\xi_E$ to
$D^2$ (lemma \ref{even}).
Thanks to the relative classifying property of $({E'},\FF_{E'})$~, this projection extends to some map
$f: H'\to E$ such that $\xi_\vp\vert H'=f^*\xi_E$~. Recall that $\xi_\vp$ and $\xi_\plat$ coincide on $\partial H'$~. Thus,
again by the relative classifying property, $\xi_\plat\vert H'$ is the pullback
of $\FF_{E'}$ by some map $g:M\to {E'}$ equal to $f$ on
 $\partial H'$~. If one likes better, rather than invoking Haefliger's classifying property \cite{haef2}, one can in this case easily build $f$ and $g$ by hands.

\bull\\

Finally, $H'$ being a handlebody, and $\pi_2(E)$, $\pi_3(E)$
being both trivial, necessarily $f$ and $g$ are homotopic rel. $\partial H'$. So $\xi_\vp$ and $\xi_\plat$ are homotopic in $H'$ rel. $\partial H'$.
This completes the homotopy argument in the twisted case, and the proof of theorem \ref{reg}.

\end{document}